# Networks of Queues Models with Several Classes of Customers and Exponential Service Times

**M. A. M. Ferreira**

Instituto Universitário de Lisboa (ISCTE – IUL), BRU - IUL, ISTAR-IUL, Lisboa, Portugal

**Abstract**

The main target of this paper is to present the Markov chain $C$ that, not giving explicitly the queue lengths stationary probabilities, has the necessary information to its determination for open networks of queues with several classes of customers and exponential service times, allowing to overcome ingeniously this problem. The situation for closed networks, in the same conditions, much easier is also presented.

## 1 Introduction

A network of queues is a collection of nodes, arbitrarily connected by arcs that are instantaneously traversed by customers, such that

 -Each node is associated with an arrivals process;

 - There is a commutation process that commands the various customers' paths.

The number of network nodes is called $J$. The arrivals processes may be composed of exogenous arrivals, from the outside of the collection, and of endogenous arrivals, from the other collection nodes. A network of queues is open if every customer may enter or leave the network. It is closed if it has a fixed number of customers that travel from node to node, there having neither arrivals from the outside of the network nor departures. The networks open from some customers and closed for others are said mixed.

The networks of queues models under study in this paper fulfil the following properties

- $J < \infty$;
- $L_j < \infty, j = 1,2, ... , J$, being $L_j$ the waiting capacity at node $j$ ;

- The customers may be of various types or, as it is also said, to belong to various classes;
- The customers of a given type, or class, arrive to the open networks according to a Poisson process;
- Every customer path in the network is fixed;
- Every service times are independent and identically and exponentially distributed.

In this situation, contrarily to what it is usual, see [5], it is not possible to approach directly the queues length process since now it is not a Markov process. Alternatively the state of a node is defined by a vector that supplies information on the position of each customer of a given class in the node and the stage in the respective path. With those assumptions, a stochastic process with the specified states is a Markov chain form which it is possible to obtain the queue lengths probabilities, see [1] and [8].

In the closed networks it is allowed that a customer changes from one class to another, when goes from one node to the following in its path. Changes during the sojourn in a node are not allowed.

## 2 Open Networks of Queues

In the open networks case it is supposed that

-There are $I < \infty$ customers' types, being the path of customer in the network defined by its type;

-A type $i$ customer follows the path $r(i,1), r(i,2), \ldots, r(i,f(i))$ where $r(i,j)$ is the $j^{th}$ node visited by the type $i$ customers and $r(i,f(i))$ the last before it abandons the network;

-The type of a customer does not change while it sojourns in the network;

-The $i$ type customers arrivals process at the network is a Poisson process with parameter $v(i), i = 1,2, \ldots, I$.

At the nodes the customers have a service requirement exponentially distributed such that

-With no loss of generality, the mean service requirement is supposed unitary;

-Every customer's service requirements, at every node, are mutually independent and independent from the arrivals processes.

If in node $j$ there are $n_j$ customers in positions $1,2, \ldots, n_j$:

-The expected service requirement for the customer in position $l$ is $\mu_j^{-1}(l), l = 1,2, \ldots, n_j$ with $\mu_j(l) > 0$ if $n_j > 0$ ;

-The service velocity given to the customer in position $l$ is denoted $\gamma_j(n_j, l)$, depends on $n_j$ and $\sum_l \gamma_j(n_j, l) = 1$;

-When a customer abandons the network, when it is in position $l$, the customers in positions $l+1, l+2, \ldots, n_j$ move to positions $l, l+1, \ldots, n_j - 1$, respectively;

- When a customer arrives at node $j$ it occupies position $l$ with probability $\delta_j(n_j + 1, l)$ and the customers that were in positions $l, l+1, \ldots, n_j$ move to positions $l+1, l+2, \ldots, n_j + 1$, respectively. Evidently $\sum_l \delta_j(n_j, l) = 1$.

Call $t_j(l)$ the customer that occupies the position $l$ at node $j$ type and $s_j(l)$ the stage of its path reached when it is at node $j$. Introduce also $c_j(l) = \left(t_j(l), s_j(l)\right)$, $c_j = \left(c_j(1), c_j(2), \ldots, c_j(n_j)\right)$ and

$$C = (c_1, c_2, \ldots, c_J) \quad (2.1).$$

If $C(t) = C$, for a given $t$, one says that the network is in sate $C$ at instant $t$. Then, with the assumptions made, $\mathbb{C} = \{C(t), t \geq 0\}$ was built like a Markov chain with a countable states space. In every node there are three kinds of possible transitions:

-Transferences between nodes, even from a node for itself. It is formalized by the operator $T_{jlr}$ meaning $T_{jlr}C$ the state of the network after a customer in position $l$ at node $j$ moves to the position $r$ in its path next node, $\left(t_j(l), s_j(l) + 1\right)$;

-Departures from the node to the outside of the network. Formalized by the operator $T_{jl.}$ meaning $T_{jl.}C$ the state of the network after a customer in position $l$ at node $j = r(i, f(i))$ moves to the outside of the network after this node;

-Arrivals at a node from the outside of the network. Formalized by the operator $T^{ir}$ meaning $T^{ir}C$ the state of the network after a $i$ type customer enters the network, coming from the exterior, and occupying the position $r$ at node $r(i, 1)$.

It is convention that $(C, l)$ means that the system state is $C$ and the customer that occupies the position $l$ in queue $j$ is the next to move. The $r$ in $(r, C^\star)$ designates the position for which the customer moves in the next node it visits, being $C^\star$ the state system after that transition.

The transition rates for the Markov chain $\mathbb{C}$ are:

$$A) q(C, T_{jlr}C) = \sum_g \sum_h q\big((C,g),(h,T_{jlr}C)\big) = \begin{cases} \sum_g \sum_h \mu_j(l)\gamma_j(n_j,l)\delta_k(n_k+1,r) \\ \sum_g \sum_h \mu_j(l)\gamma_j(n_j,l)\delta_j(n_j,r) \end{cases}$$

The double sum considers the whole $g$ and $h$ such that $T_{jlr}C = T_{jgh}C$.

$$B) q(C, T_{jl\cdot}C) = \sum_g q\big((C,g),(\cdot,T_{jg\cdot}C)\big) = \sum_g \mu_j(l)\gamma_j(n_j,l)$$

(2.2).

The sum considers the whole $g$ for which $T_{jl\cdot}C = T_{jg\cdot}C$.

$$C) q(C, T^{ir}C) = \sum_h q\big((C,\cdot),(h,T^{ih}C)\big) = \sum_h v(i)\delta_k(n_k+1,r)$$

The sum considers the whole $h$ for which $T^{ir}C = T^{ih}C$.

In **A)**, in the first case, the customer abandons node $j$ and directs to another node different from $j$. In the second case the customer abandons the node $j$ and returns to it at once. This situation is often called *feedback*.

**Theorem 2.1**

For the stationary Markov process $C$, defined by the transition rates (2.2), the equilibrium distribution is

$$\pi(C) = \prod_{j=1}^{J} \pi_j(c_j)$$

where $\pi_j(c_j) = B_j \prod_{l=1}^{n_j} \dfrac{\alpha_j\big(t_j(l),s_j(l)\big)}{\mu_j(l)}$

$$B_j^{-1} = \dfrac{\sum_{n=0}^{\infty} b_j^n}{\prod_{l=1}^{n} \mu_j(l)} \quad (2.3).$$

$$b_j = \sum_{i=1}^{I}\sum_{s=1}^{f(i)} \alpha_j(i,s)$$

and $\alpha_j(i,s) = \begin{cases} v(i), \text{if } r(i,s)=j \\ 0, \text{in other cases} \end{cases}$

**Dem.:** Begin to check that the sum of the whole $\pi(C)$ is 1, owing to the constants $B_1, B_2, \ldots, B_j$ definition. Then it is shown that the $C$ reverse process[1], $C^r$, corresponds to the entrance of $i$ type customers in the system, according to a Poisson process at rate $v(i)$, following the path $r(i,f(i)), r(i,f(i)-1), \ldots, r(i,1)$ before abandon the system. In the $C^r$ process the queues behave as in $C$ with the roles of the functions $\gamma_j$ and $\delta_j$ reciprocally changed. That is: $C^r$ has the same shape as $C$ but with different parameters. Then it is possible to define the rates corresponding to those of (2.2) for $C^r$. Finally, using the rates for $C$ and $C^r$ and $\pi(C)$ as given in (2.3), with the help of the reverse processes properties, it is possible to conclude that the proposed distribution $\pi(C)$ is the equilibrium distribution. ∎

---

[1]For the reverse process see, for instance, [1,2].

**Corollary 2.1**

The stationary queue lengths $N_j(t), j = 1,2,\ldots,J$ are, for each $t$, mutually independent and, in equilibrium, their stationary probabilities are

$$P[N_j(t) = n_j] = \pi_j(n_j) = \frac{B_j b_j^{n_j}}{\prod_{l=1}^{n_j} \mu_j(l)} \quad (2.4).$$

And, if a customer occupies the position $l$ at node $j$, the probability that it is $i$ type and is in stage $s$ is $\frac{\alpha_j(i,s)}{b_j}$. ∎

**Observation**:

-This result is a corollary of Theorem 2.1 and shows how $C$, although does not give explicitly the queue lengths stationary probabilities, has the necessary information to its determination.

-In equilibrium the nodes queue lengths are independent for each $t$: $\pi(C)$ has product form, and behave as each node was a $M|M|\cdot$ queue with arrivals rate $\alpha_j$. The queue lengths in $t_1$ and $t_2$ are independent.

# 3 Closed Networks of Queues

For the closed networks of queues model considered in this work it is assumed that

-$L_j = \infty, j = 1,2,\ldots,J$;

-The queue discipline is *FCFS*;

- The service times are independent and exponentially distributed with parameter $\mu_j$, depending on node $j$ but not on the customer type $i$.

-In equilibrium the *traffic equations*[2],

$$\alpha_k(i') = \sum \alpha_j(i'') p(j, i''; k, i'), 1 \leq k \leq J, 1 \leq i' \leq I, (k, i') \in Z_i, 1 \leq i \leq I \quad (3.1)$$

where the sum is over the whole $(j, i'') \in Z_i$ (the pairs in $Z_i$ are such that for each j, *i* only assumes the values corresponding to the customers types that can be served in node j) are fulfilled.

---

[2]The network traffic equations are easy to interpret having in account that $\alpha_j(i'')$ is the $i''$ type customers departures from node $j$ total rate and $p(j, i''; k, i')$ is the proportion of those customers that reach the node $k$ as $i'$ type customers. Summing for the whole $(j, i'') \in Z_i$, as there are no exogenous arrivals, it is obtained the $i'$ type customers arrivals total rate at node $k$, $(k, i') \in Z_i$.

The important result is:

### Theorem 3.1

If the arrivals rates satisfy the *traffic equations* (3.1), in stationary regime

$$\pi(i_1, i_2, \ldots, i_J) = B_N \prod_{j=1}^{J} \prod_{i=1}^{I} \frac{n_j!}{n_j(i)!} \left(\frac{\alpha_j(i)}{\mu_j}\right)^{n_j(i)} \quad (3.2)$$

being $B_N$ the normalizing constant, $N = \sum_{j=1}^{J} n_j$ the total number of customers in the network and $n_j(i)$ the number of $i$ type customers in node j, $i = 1,2,\ldots,I; j = 1,2,\ldots,J$.

**Dem.**: Reduces to check that in the conditions of the theorem it possible to calculate $B_N$ so that the sum of the whole probabilities in (3.2) is 1. ∎

**Observation**:

-The computation of $B_N$ may be a very hard task when the states space cardinal is very great.

## 4 Conclusions

The networks of queues with several classes of customers are very important in applications. For the open networks an example is the study of healthcare facilities. For the closed networks, an important example is the study of biological systems.

Evidently its mathematical study is much more difficult than in the case of networks with only one class of customers. The methodology presented in this work is an ingenious way to overcome the problem in the case of exponential service times, allowing to determine at least the equilibrium distribution for the population process.